\documentclass[12pt,twoside]{article}
\usepackage{amsmath,amsbsy,amsfonts}
\usepackage{hyperref}
\usepackage{color}

\newtheorem{lem}{Lemma}[section]
\newtheorem{prop}{Proposition}[section]
\newtheorem{thm}{Theorem}[section]

\newtheorem{cor}{Corollary}[section]
\newtheorem{rmk}{Remark}[section]

\newtheorem{defi}{Definition}[section]

\newcommand{\cqd}{\hspace{10pt}\fbox{}}

\let\Section=\section

\def\section{\setcounter{equation}{0}\Section}

\def\nd{\noindent}

\def\proof{{\rm \bf Proof}}

\newcommand{\w}{W_0^{1,\Phi}(\Omega)}

\begin{document}
\title{{ Quasilinear elliptic systems  with convex-concave singular  terms $\Phi\!-\!Laplacian$ operator }}

\author{
{\large Jos\'e V. Gon\c calves~~~~
\large Marcos L. Carvalho}\\
\\
{\it\small  Universidade Federal de 
Goi\'as,  Instituto de Matem\'atica e Estat\'istica}\\
{\it\small  74001-970 Goi\^ania, GO - Brazil}\\
\\
{\large Carlos Alberto Santos}\footnote{Carlos Alberto Santos acknowledges
the support of CAPES/Brazil Proc.  $N^o$ $2788/2015-02$,}\\
{\it\small  Universidade de Bras\'ilia, Departamento de Matem\'atica}\\
{\it\small   70910-900, Bras\'ilia - DF - Brazil}\\
%{\it\small e-mails:  goncalves.jva@gmail.com, csantos@unb.br,
%marcos\_leandro\_carvalho@ufg.br}\vspace{1mm}\\
}
\date{}

%\pretolerance10000

\maketitle

\begin{abstract}
\noindent {\small 
This paper deals with existence of positive solutions for a class of quasilinear elliptic systems involving the $\Phi$-Laplacian operator  and convex-concave singular terms. Our approach is based on the generalized  Galerkin Method  along with perturbartion techniques and comparison arguments in the setting of Orlicz-Sobolev spaces.}
\vskip.2cm

\end{abstract}

\section{Introduction}

This paper deals with the existence of solutions of elliptic systems of the form

\begin{equation}\label{prob}
	\left\{\
	\begin{array}{l}
		\displaystyle-\Delta_\Phi u=\frac{a_1(x)}{u^{\alpha_1}v^{\beta_1}}+b_1(x)u^{\gamma_1} v^{\sigma_1}~\mbox{in}~\Omega,
		\\
		\displaystyle-\Delta_\Phi v=\frac{a_2(x)}{u^{\beta_2}v^{\alpha_2}}+b_2(x) u^{\sigma_2} v^{\gamma_2}~\mbox{in}~\Omega,\\ 
		u,v>0~\mbox{in}~\Omega,~~u=v=0~\mbox{on}~\partial \Omega,
	\end{array}
	\right.
\end{equation}

\nd where  $\Omega\subset\mathbb{R}^N$ is a bounded domain with smooth boundary $\partial \Omega$; $\alpha_i, \beta_i, \gamma_i, \sigma_i \geq 0$ are real  constants; and  $a_i,b_i:\Omega\rightarrow \mathbb{R}$ are non-negative, mensurable functions.  
\nd In addtion, $\Phi$ is the N-function  defined by 
$$
\Phi(t)=\int_0^ts\phi(\vert s \vert)ds,~t\in\mathbb{R},
$$
\nd where $\phi:(0,\infty) \rightarrow (0,\infty)$ is   $C^1$ and satisfies: 
\begin{itemize}
	\item[($\phi_1$)]  \     \ $\mbox{(i)} \ \ s\phi(s)\to 0 \ \mbox{as} \  s\to 0,~~ \mbox{(ii)} \  s\phi(s)\to\infty \ \mbox{as} \  s\to\infty$
	\item[($\phi_2$)]  \ $s\mapsto s\phi(s) \ \mbox{\it is strictly increasing in}~ (0, \infty)$,
	\item[($\phi_3$)] there exist $\ell,m\in(1,N)$ such that 
\begin{equation*}
\ell-1\leq \frac {(s\phi(s))^\prime}{\phi(s)}\leq m-1,~s>0,
\end{equation*}
\end{itemize}
and $-\Delta_{\Phi}:\w\rightarrow W^{-1,\widetilde \Phi}(\Omega)$ is defined by 
$$
\langle -\Delta_{\Phi}u,v\rangle:=\int_\Omega \phi(|\nabla u|)\nabla u\nabla v dx,~u,v\in\w,
$$
where $\w$ stands for the classical Orlicz-Sobolev space, $W^{-1,\widetilde \Phi}(\Omega)$ denotes its dual space, and $\widetilde{\Phi}$, given by 
$$
\widetilde{\Phi}(t) = \displaystyle \max_{s \geq 0} \{ts - \Phi(s) \},~ t \geq 0,
$$ is the N-function complementary to the N-function $\Phi$ and vice-versa. We refer the reader to Section 5 for more details about this space and about the Orlicz spaces that will be denoted by  $L_{{\Phi}}(\Omega)$.

\nd In this context,  we prove an existence result of positive solutions to Problem \eqref{prob} and a weak Comparison Principle for the $\Phi-Laplacian$ operator. In this work, $\displaystyle d(x)=\inf\{ |x-y|~/~ y\in\partial \Omega\}$ for $x\in\Omega$ will stand for the distance function to the boundary of the domain $\Omega$.

\begin{thm}\label{Teor-prin}
Assume $(\phi_1)-(\phi_3)$, $0\neq a_i \in L_{\widetilde{\Psi}}(\Omega)$,  $0\neq b_i\in L^{q_i}(\Omega)$, and $a_i+b_i>0$ a.e. in $\Omega$ hold for some N-function $\Phi<\Psi << \Phi_*$ and  $q_i\geq\ell/(\ell - \sigma_i-\gamma_i-1)$, where $0<\sigma_i+ \gamma_i<\ell-1$.  Suppose in adition that $\displaystyle {a_i}{d^{-\alpha_i-\beta_i}}\in L_{\widetilde\Psi}(\Omega) $. 
  Then there exists $(u, v) \in \w\times \w$ weak solution of  \eqref{prob}   if one of the below condition is true:
  \begin{enumerate}
  \item [$(i)$] $\beta_i = 0$  $($cooperative structure$)$,
  \item [$(ii)$] $\sigma_i=0$  $($non-cooperative structure$)$,
  \item [$(iii)$] $\alpha_i=\gamma_i=0$,  and 
  $\min\{a_i,b_i\}>0$ a.e. $x \in \Omega$ $($mixed structure$)$.
  \end{enumerate}
Besides this, there exists a $C>0$ such that $u(x),v(x)\geq Cd(x)$ a.e. $x \in \Omega$.
\end{thm}

\begin{rmk}
\label{remark1}
The item $(i)$  is true if we assume  $0<\alpha_i\leq 1$ instead of $\displaystyle {a_i}{d^{-\alpha_i}}\in L_{\widetilde\Psi}(\Omega) $.
\end{rmk}

\nd  One important tool in our approuch is the comparison principle below, which is relevant by itself. Consider the problem
\begin{equation}\label{prob-aux-comp}
\left\{
\begin{array}{l}
-\Delta_\Phi u=f(x,u)~\mbox{in}~\Omega,\\
u>0~\mbox{in}~\Omega,	u=0~\mbox{on}~\partial \Omega,
\end{array}
\right.
\end{equation}
\nd where $f:\Omega\times[0,\infty)\rightarrow \mathbb{R}$ is a  Careth\'eodory function.
\begin{defi}
An	$u\in W^{1,\Phi}(\Omega)$ is a subsolution $[solution]$ $( supersolution)$  of \eqref{prob-aux-comp} if
	$$\int_\Omega\phi(|\nabla u|)\nabla u\nabla \varphi dx\leq [=](\geq) \int_\Omega f(x,u)\varphi dx,~\varphi\in \w,~\varphi\geq 0.$$
\end{defi}

\begin{thm}\label{princ-compar}{\bf (Comparison Principle}$)$
Let $\phi$ be satisfy $(\phi_1)-(\phi_3)$. Assume that 
$$
t\longmapsto\frac{f(x,t)}{t^{\ell-1}}~ \mbox{is decreasing for a.e.}~ x \in \Omega.
$$
If $u_1,u_2\in W^{1,\Phi}(\Omega)$ are  sub and supersolution of \eqref{prob-aux-comp}, respectively, such that  $\displaystyle{u_1}/{u_2}\in L^{\infty}(\Omega)$ and $u_1\leq u_2$ in $\partial \Omega$, then $u_1\leq u_2$ a.e. in $\Omega$.
\end{thm}

\nd This result extend and improve to $\Phi-Laplacian$ operator  results that are well-known to $Laplacian$ (Br\'ezis and Oswald \cite{brezis-oswald}) and $p-Laplacian$ (Diaz and Saa \cite{dsaa} and Mohammed \cite{ahmed-1}) operators. 
\vskip.5cm

\nd Below, let us do an overview about related problems to  (\ref{prob}). 
First, we point out that there is by now an extensive literature on single-equation singular problems related to (\ref{prob}) ($\phi(t)=t^{p-2}$ with $1<p<N$), that is, to problems like
\begin{equation*}\label{model prob}
  \left\{\
  \begin{array}{l}
  \displaystyle-\Delta_p  u=  \frac{a(x)}{u^{\alpha}}+b(x)u^{\gamma} ~\mbox{in}~\Omega,\\ 
  u>0~\mbox{in}~\Omega,~u=0~\mbox{on}~\partial \Omega,
  \end{array}
  \right.
\end{equation*}
where $a,b$ are appropriate potentials and $\alpha,\gamma$ are positive real constants. As it is impossible to cite all important papers that have considered this kind of problems, let us to refer the reader to works  \cite{crandall-rab-tartar, ghergu-radu,  JST,rezende-santos,lazerMcKenna-1, ahmed-1}, and their references to highlight the variety of techniques that  have been used to solve them. 

\nd Again, for the function $\phi(t)=t^{p-2}$ with $1<p<N$, we note that Problem (\ref{prob}) read as
\begin{eqnarray}
\label{sys}
\left\{\
	\begin{array}{l}
		\displaystyle-\Delta_p u=\frac{a_1(x)}{u^{\alpha_1}v^{\beta_1}}+b_1(x)u^{\gamma_1} v^{\sigma_1}~\mbox{in}~\Omega,
		\\
		\displaystyle-\Delta_p v=\frac{a_2(x)}{u^{\beta_2}v^{\alpha_2}}+b_2(x) u^{\sigma_2} v^{\gamma_2}~\mbox{in}~\Omega,\\ 
		u,v>0~\mbox{in}~\Omega,~~u=v=0~\mbox{on}~\partial \Omega,
	\end{array}
	\right.
\end{eqnarray}
and the assumptions $(\phi_1)-(\phi_3)$ are true with $\ell=m=p$. When $p=2$, problems like these have been considered with more frequency. We quote \cite{alves-correa-gonc,choi,ghergu,hai,vega,ni}, and references therein.

\nd About more general operators related to Problem (\ref{prob}), we refer  the reader to \cite{manasevich-1} and their references. An important and recent paper in this constext is \cite{JST1}. In it the authors considered a $(p,q)-Laplacian$ system and proved existence and uniqueness results for $b_i=0$ and $a_i=1$ in (\ref{sys}) by using monotonicity methods.

\nd Other classes of functions $\phi$, which satisfy $(\phi_1)-(\phi_3)$ are:
\begin{description}
\item{$\rm {(i)}$} $\phi(t)=t^{p-2}+t^{q-2}$ with $1<p<q<N$. In this case, the problem \eqref{prob}   becomes in the $(p,q)$-Laplacian problem
\begin{eqnarray*}
\left\{\
	\begin{array}{l}
		\displaystyle-\Delta_p u-\Delta_q u=\frac{a_1(x)}{u^{\alpha_1}v^{\beta_1}}+b_1(x)u^{\gamma_1} v^{\sigma_1}~\mbox{in}~\Omega,
		\\
		\displaystyle-\Delta_pv-\Delta_q v=\frac{a_2(x)}{u^{\beta_2}v^{\alpha_2}}+b_2(x) u^{\sigma_2} v^{\gamma_2}~\mbox{in}~\Omega,\\ 
		u,v>0~\mbox{in}~\Omega,~~u=v=0~\mbox{on}~\partial \Omega,
	\end{array}
	\right.
\end{eqnarray*}
\nd with $\ell=p$ and $m=q$ in the assumption $(\phi_3)$,

\item{$\rm {(ii)}$} $\phi(t)=\displaystyle \sum_{i=1}^N t^{p_i-2}$, where $1<p_1<p_2<...<p_N$, and $\frac{1}{\overline{p}}=\displaystyle \frac{1}{N}\sum_{i=1}^N\frac{1}{p_i}$ with $\overline p<N$, $\ell=p_1$, and $m=p_N$ in the assumption $(\phi_3)$. In this case, the corresponding problem is
\begin{eqnarray*}
\left\{\
	\begin{array}{l}
		\displaystyle-\sum_{i=1}^N\Delta_{p_i} u=\frac{a_1(x)}{u^{\alpha_1}v^{\beta_1}}+b_1(x)u^{\gamma_1} v^{\sigma_1}~\mbox{in}~\Omega,
		\\
		\displaystyle-\sum_{i=1}^N\Delta_{p_i} v=\frac{a_2(x)}{u^{\beta_2}v^{\alpha_2}}+b_2(x) u^{\sigma_2} v^{\gamma_2}~\mbox{in}~\Omega,\\ 
		u,v>0~\mbox{in}~\Omega,~~u=v=0~\mbox{on}~\partial \Omega,
	\end{array}
	\right.
\end{eqnarray*}
\end{description}
which is known  as an anisotropic elliptic problem in the literature.
\smallskip

\nd The proof of our Theorems are organized in three sections. In section 2, we prove Theorem \ref{princ-compar}. The principal difficulty that we found it has been to prove the convexity of an operator, which was defined inspired in one introduced by Diaz \& Saa \cite{dsaa}. Thanks to assumption 
$(\phi_3)$, we were able to show this. Also, with an accurate analysis, we removed the hypothesis $u_2/u_1 \in L^{\infty}(\Omega)$ that was considered in former results to $p-Laplacian$ operator. This was crucial in our approach.
\smallskip

\nd In section 3, we consider one ``regularization" of problem (\ref{prob}), and we proved Theorem \ref{solu-aux} that has mathematical interest by itself. In particular, we proved in this theorem that the problem
$$
\left\{\
\begin{array}{l}
-\Delta_\Phi u=\frac{{a}_1(x)}{(u + 1)^{\alpha_1}(v + 1)^{\beta_1}}+{b}_1(x)u^{\gamma_1}v^{\sigma_1},~ \mbox{in}~\Omega,\\ 
-\Delta_\Phi v=\frac{{a}_2(x)}{(u + 1)^{\beta_2}(v+ 1)^{\alpha_2}} +{b}_2(x)u^{\sigma_2}v^{\gamma_2}~ \mbox{in}~\Omega,\\  
u,v \geq 0~\mbox{in}~\Omega,~u=v=0~\mbox{on}~\partial \Omega
\end{array}
\right.
$$
admits a non-trival solution in $\w \times \w$ for $0 \leq a_i,b_i \in L^{\infty}(\Omega)$ with $a_i + b_i \neq 0$, and $\alpha_i,\beta_i, \gamma_i,\sigma_i\geq 0$ with $\gamma_i+\sigma_i < \ell -1,~i=1,2$.

\nd To prove Theorem \ref{solu-aux}, we developed an generalized Galerkin Method inspired in an idea found in the work of Browder \cite{Browder} of 1983. We believe that this approach can be very useful principally to solve problems that do not have Variational structure. In particular, we can have Variational structure in Problem (\ref{prob}) if we require a number of relationship among the powers. 
\smallskip

\nd Finally, in section 4, we complete the proof of Theorem \ref{Teor-prin}. Due our approach, we were able to consider cooperative (the right side of the first equation is increasing in $v>0$ and the right side of the second equation is increasing in $u>0$), non-cooperative and mixed structures  with a few modifications in their proofs. The most important issue in this section is to show that the sequence obtained from Theorem \ref{solu-aux} is bounded in $\w \times \w$ and it is also bounded from below by a positive function, namely, the distance function.  
\smallskip

\nd We also point out that as one motivation for us to prove the Theorem  \ref{Teor-prin}, it was the absence in literature of existence results of positive solutions to singular systems problems in Orlicz-Sobolev settings. Even in Sobolev settings for operators of the kind $p-Laplacian$ with $p\neq 2$, existence results like Theorem \ref{Teor-prin} are not frequent.

\section{Proof of Theorem \ref{princ-compar}}\label{Comp-princ}

The proof is motivated by arguments in   D\'iaz and Saa \cite{dsaa} .

\nd\proof: To begin let $J:L^1(\Omega)\rightarrow (-\infty,\infty]$ be given by 
	$$
	J(u):=
	\begin{cases}
		\displaystyle\int_{\Omega}\Phi(|\nabla u^{1/\ell}|) ,&~u\geq 0,~u^{1/\ell}\in\w,\\
		\infty, &~ \mbox{otherwise}.
	\end{cases}
	$$
We claim that the efective domain of $J$ is not empty, that is, 
$$D(J):=\{u\in L^1(\Omega)~|~u\geq 0,~u^{1/\ell}\in\w\}\neq\emptyset,$$
and, in particular, $J\not\equiv \infty$. Indeed, given $x_0\in\Omega$ take $\epsilon>0$ such that $B_\epsilon(x_0)\subset \Omega$ ($B_\epsilon(x_0)$ is the ball centered at $x_0$ and radius $\epsilon>0$) and consider the function  
$$
v_\epsilon(x):=
\begin{cases}
\overline v_\epsilon^s(|x-x_0|), & x\in B_\epsilon(x_0),\\
0,& x\in \Omega-B_\epsilon(x_0),
\end{cases}
$$ 
where $\overline v_\epsilon:[0,\epsilon]\rightarrow\mathbb{R}$ is defined by
$$
\overline v_\epsilon(t):=
\begin{cases}
1, & t=0,\\
\mbox{linear}, & 0<t<\frac\epsilon 2,\\
0& \frac\epsilon 2\leq t\leq \epsilon
\end{cases}
$$
 such that  $s>{\ell}$.

\nd So, by using Lemma \ref{lema_naru} in the Appendix, we obtain
\begin{eqnarray}
	\int_{\Omega}\Phi(|\nabla (v_\epsilon)^{1/\ell}|)dx & = & \int_A\Phi\left(\frac{s}{\ell}\overline v_\epsilon^{\frac s\ell-1}|\nabla \overline v_\epsilon|\right)dx\nonumber\\
	&\leq & C\int_A\max\left\{\overline v_\epsilon^{(\frac{s}{\ell}-1) \ell},\overline v_\epsilon^{(\frac s\ell-1)m}\right\}\Phi(|\nabla \overline v_\epsilon|)dx\nonumber\\
	&=& C\int_A \overline v_\epsilon^{(\frac s\ell-1)\ell}\Phi(|\nabla \overline v_\epsilon|)dx\leq  C\int_A  \Phi(|\nabla \overline v_\epsilon|)dx<\infty.\nonumber
\end{eqnarray}
\nd This shows our claim.

\nd Below, let us show that $J$ is a convex functional.  Set $z_i:=w_i^{1/\ell},~i=1,2$ and  $z_3=(\tau w_1+(1-\tau)w_2)^{1/\ell}$ for $w_1,w_2\in D(J)$ and $\tau \in [0,1]$. So, by applying the  H\"older inequality, we get
\begin{eqnarray}
\label{conv}
	z_3^{\ell-1}|\nabla z_3| & \leq& \tau^{\frac{\ell-1}{\ell}}z_1^{\ell-1}\tau ^{\frac{1}{\ell}}|\nabla z_1|+(1-\tau)^{\frac{\ell-1}{\ell}}z_2^{\ell-1}(1-\tau) ^{\frac{1}{\ell}}|\nabla z_2|\nonumber\\
	&\leq& (\tau z_1^\ell+(1-\tau)z_2^\ell)^{\frac{\ell-1}{\ell}}(\tau|\nabla z_1|^\ell+(1-\tau)|\nabla z_2|^\ell)^{\frac{1}{\ell}}\\
	&=& z_3^{\ell-1}(\tau|\nabla z_1|^\ell+(1-\tau)|\nabla z_2|^\ell)^{\frac{1}{\ell}}.\nonumber
\end{eqnarray}

Besides this, by computing, we get
\begin{eqnarray}
\frac{d^2}{dt^2}\Phi(t^{1/\ell})&=&\frac 1 \ell\left[(s\phi(s))^\prime{\big|_{s=t^{1/\ell}}}\frac 1 \ell(t^{1/\ell-1})^2+\phi(t^{1/\ell})\left(\frac 1\ell-1\right)(t^{1/\ell-1})^2\right]\nonumber\\
&=& \frac 1 \ell(t^{1/\ell-1})^2\phi(t^{1/\ell})\left[\frac{(s\phi(s))^\prime}{\phi(s)}\Big|_{s=t^{1/\ell}}\frac 1\ell+\frac 1\ell-1\right]\nonumber\\
&\stackrel{(\phi_3)}\geq& \frac 1 \ell(t^{1/\ell-1})^2\phi(t^{1/\ell})\left[(\ell-1)\frac 1\ell+\frac 1\ell-1\right]=0,\nonumber
\end{eqnarray}
because we used the hypothesis $(\phi_3)$ to obtain the last inequality. That is,  $\Phi(t^{1/\ell})$ for $t>0$ is a convex function. 

\nd So, it follows from (\ref{conv})  and  of the convexity of $t\mapsto \Phi(t^{\frac{1}{\ell}})$ for $t>0$,  that
\begin{eqnarray*}
	J(\tau w_1+(1-\tau)w_2) &=&\int_{\Omega}\Phi(|\nabla z_3|)dx\nonumber\\
	&\leq & \int_{\Omega}\Phi\left((\tau|\nabla z_1|^\ell+(1-\tau)|\nabla z_2|^\ell\right)^{\frac{1}{\ell}})dx\nonumber\\
	&\leq & \tau J(w_1)+(1-\tau)J(w_2),
\end{eqnarray*}
showing that $J$ is a convex functional. 

\nd Now, if we assumed that 
$
\Omega_0:=\{x\in\Omega~|~u_1(x)>u_2(x)\}
$
has positive Lebesgue meassure, then  $\displaystyle\varphi_i=({u_1^{\ell}-u_2^\ell})^+/{u_i^{\ell-1}},~i=1,2$ would be non-null admissible test functions, because
$
\varphi_1=\varphi_2=0~\mbox{in}~ \Omega_0^c,
$
$$
\vert \nabla \varphi_i\vert \leq \ell \Vert {u_j}/{u_i} \Vert^{\ell-1}_{\infty}\vert \nabla \varphi_j\vert + [1 + (\ell -1)\Vert {u_j}/{u_i} \Vert^{\ell}_{\infty}]\vert \nabla \varphi_i\vert ,
$$
 ${u_i}/{u_j}\in L^{\infty}(\Omega_0)$ for $i\neq j$, Lemma \ref{lema_naru}, and convexity of $\Phi$.

\nd So, it follows from the convexity of $J$, by using $\varphi_i$ as test functions and the fact that $u_1,u_2\in D(J)$, that
\begin{eqnarray}\label{diaz2}
0&\leq&\langle J^\prime(u_1^\ell)-J^\prime(u_2^{\ell}),u_1^\ell-u_2^\ell\rangle\nonumber\\
&=&\int_\Omega\phi(|\nabla u_1|)\nabla u_1\nabla \left(\frac{u_1^\ell-u_2^\ell}{u_1^{\ell-1}}\right)-\phi(|\nabla u_2|)\nabla u_2\nabla \left(\frac{u_1^\ell-u_2^\ell}{u_2^{\ell-1}}\right)dx.
\end{eqnarray}

\nd Since $u_1$ is a subsolution and $u_2$ is a supersolution of problem (\ref{prob-aux-comp}), it follows from  \eqref{diaz2} and the fact that $\displaystyle t\mapsto{f(x,t)}/{t^{\ell-1}}$ is decreasing that
\begin{eqnarray*}
0 &\leq& \int_{\Omega}[\phi(|\nabla u_1|)\nabla u_1\nabla \varphi_1-\phi(|\nabla u_2|)\nabla u_2\nabla \varphi_2]dx\nonumber\\
&=&\int_{\Omega}\left[\phi(|\nabla u_1|)\nabla u_1\nabla \left(\frac{u_1^\ell-u_2^\ell}{u_1^{\ell-1}}\right)-\phi(|\nabla u_2|)\nabla u_2\nabla \left(\frac{u_1^\ell-u_2^\ell}{u_2^{\ell-1}}\right)\right]dx\nonumber\\
&\leq&\int_{\Omega_0}\left(f(x,u_1)\varphi_1-f(x,u_2)\varphi_2\right)dx\nonumber\\
&=&\int_{\Omega_0}\left(\frac{f(x,u_1)}{u^{\ell-1}}-\frac{f(x,u_2)}{u^{\ell-1}}\right)(u_1^{\ell}-u_2^{\ell})dx<0,
\end{eqnarray*}
but this is impossible, that is, $\Omega_0$ has null Lebesgue meassure. This ends our proof. \hfill\cqd

\section{Problem (\ref{prob})  regularized}

Let us regularize Problem (\ref{prob}) by summing $\varepsilon>0$ in singular term, and adding  $\delta \geq 0$ on non-sigular term. The last one is to easy the application of Theorem \ref{princ-compar} in the proof of Theorem \ref{Teor-prin}.

\nd Let $\epsilon \in (0,1)$ and $\delta\geq 0$. So, we associate with (\ref{prob}) the  ``regularized" problem 
\begin{equation}\label{auxprob}
\left\{\
\begin{array}{l}
-\Delta_\Phi u=\frac{\hat{a}_1(x)}{(u + \epsilon)^{\alpha_1}(v + \epsilon)^{\beta_1}}+\hat{b}_1(x)(u+\delta)^{\gamma_1}(v+\delta)^{\sigma_1},~ \mbox{in}~\Omega,\\ 
-\Delta_\Phi v=\frac{\hat{a}_2(x)}{(u + \epsilon)^{\beta_2}(v|+ \epsilon)^{\alpha_2}} +\hat{b}_2(x)(u+\delta)^{\sigma_2}(v+\delta)^{\gamma_2}~ \mbox{in}~\Omega,\\  
u,v \geq 0~\mbox{in}~\Omega,~u=v=0~\mbox{on}~\partial \Omega.
\end{array}
\right.
\end{equation}

\nd So, we have.
\begin{thm}\label{solu-aux}
	Assume $(\phi_1)-(\phi_3)$ and $0 \leq \hat{a}_i,\hat{b}_i\in L^{\infty}(\Omega)$ hold with $\hat{a}_i + \hat{b}_i\neq 0$. Suppose in adition that  $\alpha_i,\beta_i,\gamma_i,\sigma_i\geq 0$ with $ \gamma_i+\sigma_i <\ell -1$.
	Then there exists a weak solution $(u^1,u^2)=(u^1_{\varepsilon,\delta},u^2_{\varepsilon,\delta})\in \w\times \w$ of Problem  \eqref{auxprob} for each $\varepsilon>0$ and $\delta\geq 0$ given. Beside this, $u^1,u^2\neq 0$.
\end{thm}
\proof:  Consider the vector space 
$
E:=\w\times\w
$
\nd endowed with the norm 
$
\|(u,v)\|:=\|u\|+\|v\|.
$
So, $(E,\|.\|)$ is a reflexive Banach space. (We refer the reader to Section  5 for for some basic facts on Orlicz-Sobolev spaces as well as references). 

\nd Consider the mapping $A=A_{\varepsilon,\delta}:E\times E\longrightarrow\mathbb{R}$ defined by
\begin{eqnarray}\label{operadorA}
A(u_1,u_2,\varphi,\psi)& = & \int_\Omega \big[ \phi(|\nabla u_1|)\nabla u_1\nabla \varphi+ \phi(|\nabla u_2|)\nabla u_2\nabla \psi\big ]dx\nonumber\\
& - & \int_\Omega \frac{\hat{a}_1(x)\varphi}{(|u_1| + \epsilon)^{\alpha_1}(|u_2| + \epsilon)^{\beta_1}}+\frac{\hat{a}_2(x)\psi}{(|u_1|+\epsilon)^{\beta_2}(|u_2|+\epsilon)^{\alpha_2}}dx\nonumber\\
& - & \int_\Omega \hat{b}_1(x)(u_1^++\delta)^{\gamma_1}(u_2^++\delta)^{\sigma_1} \varphi dx\\
& - &\int_\Omega \hat{b}_2(x)(u_1^++\delta)^{\sigma_2}(u_2^++\delta)^{\gamma_2} \psi dx.\nonumber
\end{eqnarray}

\nd So, we have.

\begin{prop}\label{A-bem-defi}
	The functional $A$  in \eqref{operadorA} is well-defined, linear  and satisfies:  
$$
A(u_1,u_2,.,.)\in (E\times E)'=E' \times E'~ \mbox{for each}~ (u_1,u_2)\in E.
$$	
\end{prop}
\proof:  Given  $(u_1,u_2) \in E$,  let $A_1=A_{1,\varepsilon,\delta}$ be given by
\begin{eqnarray}
	A_1(u_1,u_2,\varphi)&:=&\int_{\Omega} \phi(|\nabla u_1|)\nabla u_1\nabla \varphi-\frac{\hat{a}_1(x)\varphi}{(|u_1|+\epsilon)^{\alpha_1}(|u_2|+\epsilon)^{\beta_1}}dx\nonumber\\
	&-& \int_{\Omega}\hat{b}_1(x)(u_1^++\delta)^{\gamma_1}(u_2^++\delta)^{\sigma_1}\varphi dx,\nonumber
\end{eqnarray}
for $\varphi\in \w$ and $A_2=A_{2,\varepsilon,\delta}$ defined by
\begin{eqnarray}
	A_2(u_1,u_2,\psi)&:=&\int_{\Omega} \big[ \phi(|\nabla u_2|)\nabla u_2\nabla \psi - \frac{\hat{a}_2(x)\psi}{(|u_1|+\epsilon)^{\beta_2}(|u_2|+\epsilon)^{\alpha_2}}dx\nonumber\\	
	&-&\int_\Omega \hat{b}_2(x)(u_1^++\delta)^{\sigma_2}(u_2^++\delta)^{\gamma_2} \psi dx,\nonumber
\end{eqnarray}
for $\psi\in \w$. That is, we are rewriting $A$ as  $A:=(A_1,A_2)$. \\

\nd {\bf Claim:} $A_1(u_1,u_2,\varphi)$ is well-defined. It follows from  H\"older's inequality, the embedding $\w\hookrightarrow L_\Phi(\Omega)$, the inequality $\widetilde{\Phi}(t\phi(t))\leq {\Phi}(2t)$ and the fact that $\Phi\in\Delta_2$, that
\begin{eqnarray}\label{A-bem-def}
	|A_1(u_1,u_2,\varphi)|&\leq&\int_{\Omega}\phi(|\nabla u_1|)|\nabla u_1||\nabla \varphi|+\frac{\hat{a}_1(x)|\varphi|}{\epsilon^{\alpha_1+\beta_1}} dx\nonumber\\
	& + & \int_\Omega \hat{b}_1(x)(u_1^++\delta)^{\gamma_1}(u_2^++\delta)^{\sigma_1}|\varphi|dx\nonumber\\
	&\leq& 2\|\phi(|\nabla u_1|)|\nabla u_1|\|_{\widetilde \Phi}\| \varphi\|+\frac {2}{\epsilon^{\alpha_1+\beta_1}}\|\hat{a}_1\|_{\infty}\|\varphi\|_\Phi\nonumber\\
	& + & 2\|\hat{b}_1(|u_1|+\delta)^{\gamma_1}(|u_2|+\delta)^{\sigma_1}\|_{\widetilde \Phi}\|\varphi\|_\Phi\nonumber\\
	&\leq& 2\|\phi(|\nabla u_1|)|\nabla u_1|\|_{\widetilde \Phi}\| \varphi\|+\frac {C}{\epsilon^{\alpha_1+\beta_1}}\|\hat{a}_1\|_{\infty}\| \varphi\|\nonumber\\
	& + & C\|\hat{b}_1(|u_1|+\delta)^{\gamma_1}(|u_2|+\delta)^{\sigma_1}\|_{\widetilde \Phi}\| \varphi\|.\nonumber
\end{eqnarray}
To end the proof of the Claim, it remains to show that  
$$\|\hat{b}_1(|u_1|+\delta)^{\gamma_1}(|u_2|+\delta)^{\sigma_1}\|_{\widetilde \Phi}<\infty.$$
 Indeed, since $L_\Phi(\Omega)\hookrightarrow L^{{(\gamma_1+\sigma_1)\ell}/({\ell-1})}(\Omega)$, because $L_\Phi(\Omega)\hookrightarrow L^{\ell}(\Omega)$ and ${\gamma_1+\sigma_1}\in(0,\ell-1)$,  we obtain by Lemma \ref{lema_naru_*}, that
\begin{eqnarray}
\label{551}
\displaystyle\int_\Omega \widetilde\Phi(\hat{b}_1(x)(|u_1|\!\!&\!\!\!\!\!\!\!+\!\!\!\!&\!\!\!\!\!\delta)^{\gamma_1}(|u_2|+\delta)^{\sigma_1}) \nonumber\\
	&\leq& \max\{\|\hat{b}_1\|_\infty^{\frac{\ell}{\ell-1}},\|\hat{b}_1\|_\infty^{\frac{m}{m-1}}\}\int_\Omega \widetilde\Phi((|u_1|+|u_2|+\delta+1)^{\gamma_1+\sigma_1})dx\nonumber\\
	&\leq& C\int_\Omega(|u_1|+|u_2|+\delta+1)^{\frac{(\gamma_1+\sigma_1)\ell}{\ell-1}}dx \\
	&\leq & C \Vert |u_1|+|u_2|+\delta+1\Vert_{\Phi}^{\frac{(\gamma_1+\sigma_1)\ell}{\ell-1}} <\infty	,\nonumber
\end{eqnarray}
that is, 
\begin{equation}
\label{56}
|A_1(u_1,u_2,\varphi)|\leq [2\|\phi(|\nabla u_1|)|\nabla u_1|\|_{\widetilde \Phi}+\frac {C_1}{\epsilon^{\alpha_1+\beta_1}}\|\hat{a}_1\|_{\infty}+ C_2\Vert(u_1,u_2)\Vert^{\frac{(\gamma_1+\sigma_1)\ell}{\ell-1}}+C_\delta]\| \varphi\| .
\end{equation}
 In a similar way one shows that  
\begin{equation}
\label{57}
|A_2(u_1,u_2,\psi)|\leq [2\|\phi(|\nabla u_2|)|\nabla u_2|\|_{\widetilde \Phi}+\frac {D_1}{\epsilon^{\alpha_2+\beta_2}}\|\hat{a}_2\|_{\infty}+ D_2\Vert(u_1,u_2)\Vert^{\frac{(\gamma_2+\sigma_2)\ell}{\ell-1}}+D_\delta]\| \psi\| ,
\end{equation}
\nd for each $\varepsilon>0$ and $\delta\geq 0$ given, where  $C_i,D_i>0$ and $C_\delta,D_\delta\geq0$ with $C_\delta=D_\delta=0$ if $\delta=0$. The linearity is clear. These end the proof.\hfill\cqd

\begin{prop}\label{operador-T}
	There is an only operator $T=T_{\varepsilon,\delta}:E\longrightarrow E'$ such that $\langle T(u_1,u_2),(\varphi,\psi)\rangle=A(u_1,u_2,\varphi,\psi)$ for all $(u_1,u_2),(\varphi,\psi)\in E$.
\end{prop}
\proof: Let $(u_1,u_2)\in E$. Of course, it there is at most one such $T$. In addition, by \eqref{56} and \eqref{57}, we have
\begin{eqnarray*}
\|T(u_1,u_2)\|_{E'}&\leq& 2\|\phi(|\nabla u_1|)|\nabla u_1|\|_{\widetilde \Phi}\| \varphi\|+\frac {C_1}{\epsilon^{\alpha_1+\beta_1}}\|\hat{a}_1\|_{\infty}\| \varphi\|\nonumber\\
& + & C_2\Vert(u_1,u_2)\Vert^{\frac{(\gamma_1+\sigma_1)\ell}{\ell-1}}\| \varphi\|+C_\delta\| \varphi\|.
\nonumber\\
&+&  2\|\phi(|\nabla u_2|)|\nabla u_2|\|_{\widetilde \Phi}\| \psi\|+\frac {D_1}{\epsilon^{\alpha_2+\beta_2}}\|\hat{a}_2\|_{\infty}\| \psi\|\nonumber\\
&+& D_2\Vert(u_1,u_2)\Vert^{\frac{(\gamma_2+\sigma_2)\ell}{\ell-1}}\| \varphi\|+D_\delta\| \psi\|,
\end{eqnarray*}
showing that  $T(u_1,u_2)\in E'$.\hfill\cqd
\smallskip

 \nd Our next aim is to show that there exist $(u^{1},u^{2})\in E\setminus \{0\}$ such that  $T(u^{1},u^{2})=0$. In fact, we will have that this  $(u^{1},u^{2})$ will be  a non-negative weak solution of the system \eqref{auxprob}.
 
\nd Now, since $a_i+b_i\neq 0$, we can take $(\omega_1,\omega
_2)\in E$  such that 
\begin{equation}\label{w}
(\hat{a}_i+ \hat{b}_i)\omega_i\neq 0~\mbox{and}~  (\hat{a}_i+ \hat{b}_i)\omega_i\in L^1(\Omega),~i=1,2.
\end{equation}

\nd From now on, let us consider the below set of linear subspaces of $\w$, that is, 
\begin{equation}
\label{defin}
\mathcal{A}=\left\{F\subset\w~\big|~F\mbox{ is a linear subspace; }\omega_i\in F~\mbox{and}~\dim F<\infty\right\}.
\end{equation}
preordered by set inclusion.

\nd Take a such $F\in \mathcal{A}$. Let   $\beta=\{e_1,e_2,...,e_s\}$ be a linear basis of $F$, where $s:=\dim F$ is denoting the dimension of $F$. So, there exist an unique  $\xi^i=(\xi_1^i,\xi_2^i,...,\xi_s^i)$ such that 
$$
(u,v)=\sum_{j=1}^s(\xi_j^1e_j,\xi_j^2e_j),
$$
for each $u,v\in F$ given.

\nd Consider the isometric embedding
$$
\begin{array}{c}
I_F:(F\times F,\|.\|)\longrightarrow(E,\|.\|)~\mbox{defined by}~
I_F(u_1,u_2)=(u_1,u_2),
\end{array}
$$
and
$$ 
T_F:=I_F^\prime\circ T\circ I_F:F\times F\longrightarrow F'\times F',
$$ 
\nd where  $I_F^{\prime}$ is the adjoint of $I_F$. So, we have
\begin{eqnarray}
\langle T_F (u_1,u_2),(\psi,\varphi)\rangle&=&\langle I_F^\prime\circ T\circ I_F (u_1,u_2),(\psi,\varphi)\rangle=\langle T\circ I_F(u_1,u_2),I_F(\psi,\varphi)\rangle\nonumber\\
&=&\langle T (u_1,u_2),(\psi,\varphi)\rangle~\mbox{for all}~( u_1,u_2),(\psi,\varphi)\in F\times F.\nonumber
\end{eqnarray}
\nd Thus
\begin{eqnarray}\label{T_K}
	\langle T_F(u_1,u_2),(\varphi,\psi)\rangle & = & \int_\Omega \big[ \phi(|\nabla u_1|)\nabla u_1\nabla \varphi+ \phi(|\nabla u_2|)\nabla u_2\nabla \psi\big ]dx\nonumber\\
	& - & \int_\Omega \frac{\hat{a}_1(x)\varphi}{(|u_1| + \epsilon)^{\alpha_1}(|u_2| + \epsilon)^{\beta_1}}+\frac{\hat{a}_2(x)\psi}{(|u_1|+\epsilon)^{\beta_2}(|u_2|+\epsilon)^{\alpha_2}}dx\nonumber\\
	& - & \int_\Omega \hat{b}_1(x)(u_1^++\delta)^{\gamma_1}(u_2^++\delta)^{\sigma_1} \varphi dx\\
	& - &\int_\Omega \hat{b}_2(x)(u_1^++\delta)^{\sigma_2}(u_2^++\delta)^{\gamma_2} \psi dx,~\mbox{for all}~ (u_1,u_2),(\psi,\varphi)\in F\times F\nonumber
\end{eqnarray}

\nd In this context, we have.
\begin{prop}\label{S_K-cont}
	The operator $T_F$ is continuous.
\end{prop}
\proof:  To this end, we
set $T_F=(T_1,T_2)$, where
\begin{eqnarray}
	\langle T_1(u_1,u_2),\varphi \rangle &:=&\int_{\Omega} \phi(|\nabla u_1|)\nabla u_1\nabla \varphi-\frac{\hat{a}_1(x)\varphi}{(|u_1|+\epsilon)^{\alpha_1}(|u_2|+\epsilon)^{\beta_1}}dx\nonumber\\
	& - &\int_{\Omega}\hat{b}_1(x)(u_1^++\delta)^{\gamma_1}(u_2^++\delta)^{\sigma_1}\varphi dx,\nonumber\\
	\langle T_2(u_1,u_2),\psi \rangle &:= &	\int_{\Omega} \big[ \phi(|\nabla u_2|)\nabla u_2\nabla \psi - \frac{\hat{a}_2(x)\psi}{(|u_1|+\epsilon)^{\beta_2}(|u_2|+\epsilon)^{\alpha_2}}dx\nonumber\\	
	&-&\int_\Omega \hat{b}_2(x)(u_1^++\delta)^{\sigma_2}(u_2^++\delta)^{\gamma_2} \psi dx,\nonumber
\end{eqnarray}
\nd for all $(u_1,u_2) \in F\times F$ and $\varphi,\psi\in F$,
and  we note that the operator  $-\Delta_{\Phi}$ is continuous (see \cite[Lemma 3.1]{Fukagai}). Therefore, it just remains to show that 
$$
T_i-(-\Delta_{\Phi}){\big|}_F,~i=1,2
$$ 
\nd are continuous. 

\nd To show these, let 
$(u_{1,n},u_{2,n}) \subseteq F\times F$ such that $(u_{1,n},u_{2,n})\rightarrow (u_1,u_2)$ in $F\times F$.
So,  passing to a subsequence if necessary, using Lemma \ref{lema_Phi} and the embedding $L_\Phi(\Omega)\hookrightarrow L^\ell(\Omega)$, we have
\begin{itemize}
	\item[(1)] $u_{i,n}\rightarrow u_i$ a.e. in $\Omega,~i=1,2$;
	\item[(2)] there exist $h_i\in L^\ell(\Omega)$ such that $|u_{i,n}|\leq h_i,~i=1,2$. 
\end{itemize}

\nd So,
\begin{eqnarray}
	\frac{\hat{a}_1(x)\varphi}{(|u_{1,n}|+\epsilon)^{\alpha_1}(|u_{2,n}|+\epsilon)^{\beta_1}} &\stackrel{\mbox{a.e.}}\longrightarrow&\frac{\hat{a}_1(x)\varphi}{(|u_1|+\epsilon)^{\alpha_1}(|u_{2}|+\epsilon)^{\beta_1}},\nonumber\\	
	\hat{b}_1(x)(u_{1,n}^++\delta)^{\gamma_1}(u_{2,n}^++\delta)^{\sigma_1} \varphi 
	&\stackrel{\mbox{a.e.}}\longrightarrow& \hat{b}_1(x)(u_1^++\delta)^{\gamma_1}(u_2^++\delta)^{\sigma_1}  \varphi,\nonumber
\end{eqnarray}

\nd and, by using the facts that $\widetilde{\Phi}$ is convex and it satisfies $\Delta_2$, we get
\begin{eqnarray}\label{a}
   \widetilde{\Phi}\left(\left|\frac{\hat{a}_1(x)}{(|u_{1,n}|+\epsilon)^{\alpha_1}(|u_{2,n}|+\epsilon)^{\beta_1}}-\frac{\hat{a}_1(x)}{(|u_1|+\epsilon)^{\alpha_1}(|u_{2}|+\epsilon)^{\beta_1}}\right|\right) \leq  \widetilde{\Phi} \left(\frac{2|\hat{a}_1|_\infty}{\epsilon^{\alpha_1+\beta_1}}\right).\nonumber%\in L^1(\Omega).
\end{eqnarray}

\nd So, by Lebesgue's Theorem, we have
$$
\int_\Omega\widetilde{\Phi}\left(\left|\frac{\hat{a}_1(x)}{(|u_{1,n}|+\epsilon)^{\alpha_1}(|u_{2,n}|+\epsilon)^{\beta_1}}-\frac{\hat{a}_1(x)}{(|u_1|+\epsilon)^{\alpha_1}(|u_{2}|+\epsilon)^{\beta_1}}\right|\right)dx\rightarrow 0
$$
or, in an equivalent way,
\begin{equation}\label{conv-u1n}
	\left\|\frac{\hat{a}_1}{(|u_{1,n}|+\epsilon)^{\alpha_1}(|u_{2,n}|+\epsilon)^{\beta_1}}-\frac{\hat{a}_1}{(|u_1|+\epsilon)^{\alpha_1}(|u_{2}|+\epsilon)^{\beta_1}}\right\|_{\widetilde{\Phi}}\rightarrow 0,
\end{equation}
because  $\widetilde{\Phi}\in\Delta_2$. That is, by the H\"older inequality and \eqref{conv-u1n}, we obtain
$$ 
\int_\Omega\left(\frac{\hat{a}_1(x)}{(|u_{1,n}|+\epsilon)^{\alpha_1}(|u_{2,n}|+\epsilon)^{\beta_1}}-\frac{\hat{a}_1(x)}{(|u_1|+\epsilon)^{\alpha_1}(|u_{2}|+\epsilon)^{\beta_1}}\right)vdx\rightarrow0,
$$
\nd for each  $v\in\w$.  

\nd By similar arguments to the above ones, we have\\ \\
$\widetilde{\Phi}\left(\hat{b}_1(x)|(u_{1,n}^++\delta)^{\gamma_1}(u_{2,n}^++\delta)^{\sigma_1}-(u_{1}^++\delta)^{\gamma_1}(u_{2}^++\delta)^{\sigma_1}| 
\right)$
\begin{eqnarray}
	&\leq&\widetilde{\Phi}\left(2|\hat{b}_1|_\infty\frac{(h_1+\delta)^{\gamma_1}(h_2+\delta)^{\sigma_1}+(u_{1}^++\delta)^{\gamma_1}(u_{2}^++\delta)^{\sigma_1}}{2}\right)\nonumber\\	
	&\leq& C\left(\widetilde{\Phi}((h_1+h_2+\delta)^{\gamma_1+\sigma_1})+\widetilde{\Phi}((u_{1}^+ + u_{2}^+  + \delta)^{\gamma_1+\sigma_1})\right)\nonumber\\
	&\leq&	C\left((h_1+h_2+\delta)^{(\gamma_1+\sigma_1)\ell/(\ell -1)}  +  (u_{1}^+ + u_{2}^+  + \delta)^{(\gamma_1+\sigma_1)\ell/(\ell-1)}  +C_\delta\right)\in L^1(\Omega),	\nonumber
\end{eqnarray}
\nd where  $C_\delta\geq 0$. So,   we have
$$
\int_\Omega \hat{b}_1(x)[(u_{1,n}^++\delta)^{\gamma_1}(u_{2,n}^++\delta)^{\sigma_1}-(u_{1}^++\delta)^{\gamma_1}(u_{2}^++\delta)^{\sigma_1}]vdx\longrightarrow0.
$$
\nd These show that $T_1-(-\Delta_{\Phi}){\big|}_F$ is continuous. In a similar way one shows that  $T_2-(-\Delta_{\Phi}){\big|}_F$ is continuous as well. Therefore  $T_F$ is continuous, ending the proof.\hfill\cqd 
\medskip

\nd We are going to use the proposition below, which is a consequence of Brouwer's Fixed Point Theorem, see e.g.  Lions \cite{lions}.
\begin{prop}\label{prop-principal}
	Suppose that $S:\mathbb{R}^m\rightarrow\mathbb{R}^m$ is a continuous function such that $\langle S(\eta),\eta\rangle>0$ on $|\eta|=r$, where $\langle \cdot,\cdot\rangle$ is the usual inner product in $\mathbb{R}^m$ and $|\cdot|$ is its corresponding norm. Then, there exists $\eta_0\in \overline B_r(0)$ such that $S(\eta_0)=0$.
\end{prop}

\nd To apply this preposition, we have to reduce our operator $T_F$ to a finite dimensional space. To do this, define
$	S_F:=i^\prime\circ T_F\circ i:\mathbb{R}^s\times\mathbb{R}^s\rightarrow\mathbb{R}^s\times\mathbb{R}^s,$
where
 $i=i_F:(\mathbb{R}^s\times\mathbb{R}^s,|\cdot|)\rightarrow(F\times F,\|\cdot\|)$,  given by $i(\xi^1,\xi^2)=(u,v)$, is an  isometry with the norm  in $\mathbb{R}^s\times\mathbb{R}^s$ given by $|(\xi^1,\xi^2)|:=\|(u,v)\|:=\|u\|+\|v\|$, and  $i^{\prime}$ is its adjoint operator.

\begin{prop}\label{raizT_K}
	The operator  $S_F$  admits one zero $(\xi^1,\xi^2)=(\xi^{1,\varepsilon,\delta}_{F},\xi^{2,\varepsilon,\delta}_{F})$ in $\mathbb{R}^s\times\mathbb{R}^s$ with $\xi^1,\xi^2\neq 0$. Moreover, the corresponding vector $(u_1,u_2)=(u^{1,\varepsilon,\delta}_F,u^{2,\varepsilon,\delta}_F)=i(\xi^{1},\xi^{2})\in F\times F$  satisfies $T_F(u_1,u_2)=0$ and $u_1,u_2\neq 0$.
\end{prop}
\nd\proof. Given $(\xi^1,\xi^2)\in\mathbb{R}^s\times \mathbb{R}^s$, denote the by $(u_1,u_2)\in F\times F$ the only image vector of $(\xi^1,\xi^2)$ by the isometry $i$. So, by using  $(\phi_3)$, we have
\begin{eqnarray}
\label{511}
	(S_F(\xi^1,\xi^2),(\xi^1,\xi^2))&=&(i^\star\circ T_F\circ i(\xi^1,\xi^2),(\xi^1,\xi^2))
	=(T_F(u_1,u_2),(u_1,u_2))\nonumber\\
	&\geq& \int_\Omega \big[ \phi(|\nabla u_1|)|\nabla u_1|^2+ \phi(|\nabla u_2|)|\nabla u_2|^2\big ]dx\nonumber\\
	& - & \int_\Omega \frac{\hat{a}_1(x)}{\epsilon^{\alpha_1+\beta_1}}|u_1|+\frac{\hat{a}_2(x)}{\epsilon^{\alpha_2+\beta_2}}|u_2|dx\\
	& - & \int_\Omega \hat{b}_1(x)(u_1^++\delta)^{\gamma_1+1}(u_2^++\delta)^{\sigma_1}dx\nonumber\\
	& - &\int_\Omega \hat{b}_2(x)(u_1^++1)^{\sigma_2}(u_2^++1)^{\gamma_2+1} dx.\nonumber
\end{eqnarray}
Now,  by using  $(\Phi_3)$ and following similar arguments as done in (\ref{551}), we obtain
\begin{eqnarray}
	(S_F(\xi^1,\xi^2),(\xi^1,\xi^2))&\geq&\ell\min\{\|u_1\|^\ell,\|u_1\|^m\}+\ell\min\{\|u_2\|^\ell,\|u_2\|^m\}\nonumber\\
	&-& C^\varepsilon_1\|(u_1,u_2)\|-C^\delta_2\|(u_1,u_2)\|^{\gamma_1+\sigma_1+1}\nonumber\\
	&-& C^\delta_3\|(u_1,u_2)\|^{\gamma_2+\sigma_2+1} -C_4^\delta,\nonumber
\end{eqnarray}
\nd where    $C^\varepsilon_1,C^\delta_2,C^\delta_3>0$, and $C_4^\delta\geq0$ are real constants. 

\nd So, by setting $r_0:= \|(u_1,u_2)\| = r_1+r_2 \geq 2$,  we have $r_1:=\|u_1\|\geq 1$ or $r_2:=\|u_1\|\geq 1$, and
$$r_0^\ell=(r_1+r_2)^\ell\leq 2^\ell\min\{r_1^\ell,r^m_1\}+2^\ell\min\{r_2^\ell,r^m_2\},$$
that is,
\begin{eqnarray}
(S_F(\xi_1,\xi_2),(\xi_1,\xi_2)) 	&\geq& \frac{\ell}{2^\ell}r_0^\ell- C^\varepsilon_1r_0-C^\delta_2r_0^{\gamma_1+\sigma_1+1}-C^\delta_3r_0^{\gamma_2+\sigma_2+1}-C_4^\delta.\nonumber 
\end{eqnarray}
Since,  $\gamma_j+\sigma_j\in(0,\ell-1)$, we can choose an $r_0=r_0^{\varepsilon,\delta}>0$ such that 
$$(S_F(\xi_1,\xi_2),(\xi_1,\xi_2)) >0~\mbox{for all}~ \Vert (\xi_1,\xi_2) \Vert = r_0$$
which implies, by applying Proposition \ref{prop-principal}, that there is a $(\xi^1,\xi^2)=(\xi_F^{1,\varepsilon,\delta},\xi_F^{2,\varepsilon,\delta})\in \overline B_{r_0}(0,0)\subset\mathbb{R}^s\times\mathbb{R}^s$ such that $S_F(\xi^1,\xi^2)=(0,0)$. Let $(u_1,u_2)=(u_F^{1,\varepsilon,\delta},u_F^{2,\varepsilon,\delta})=i(\xi_F^{1,\varepsilon,\delta},\xi_F^{2,\varepsilon,\delta})$. Then
$$
\langle T_{F}(u_1,u_2),(\psi,\varphi)\rangle=( S_{F}(\xi^1,\xi^2),(\eta_1,\eta_2))=0,
$$
\nd for all  $(\psi,\varphi)=i(\eta_1,\eta_2) \in F \times F$, that is,
\begin{eqnarray}\label{solaux-1}
	\int_{\Omega} \phi(|\nabla u_1|)\nabla u_1\nabla \psi & = &  \int_{\Omega} \frac{\hat{a}_1(x)\psi}{(|u_1|+\epsilon)^{\alpha_1}(|u_2|+\epsilon)^{\beta_1}}dx\nonumber\\
	& + &\int_{\Omega}\hat{b}_1(x)(u_1^++\delta)^{\gamma_1}(u_2^++\delta)^{\sigma_1}\psi dx,~\psi\in F\nonumber
\end{eqnarray}
and
\begin{eqnarray}\label{solaux-2}
	\int_{\Omega} \big[ \phi(|\nabla u_2|)\nabla u_2\nabla \varphi 
	& = &  \int_{\Omega}\frac{\hat{a}_2(x)\varphi}{(|u_1|+\epsilon)^{\beta_2}(|u_2|+\epsilon)^{\alpha_2}}dx\nonumber\\	
	& + & \int_\Omega \hat{b}_2(x)(u_1^++\delta)^{\sigma_2}(u_2^++\delta)^{\gamma_2} \varphi dx,~\varphi\in F.\nonumber
\end{eqnarray}
\nd In particular, we have that   $u_1\not\equiv0$, because otherwise
\begin{eqnarray*}
	0
	& = &  \int_{\Omega}\frac{\hat{a}_1(x)\varphi}{\epsilon^{\beta_1}(|u_2|+\epsilon)^{\alpha_1}}dx+\int_\Omega \hat{b}_1(x)\delta^{\sigma_1}(u_2^++\delta)^{\gamma_1} \varphi dx,~\varphi\in F.\nonumber
\end{eqnarray*}
\nd contradicting \eqref{w}. In  a similar way, we show that $u_2\neq0$ as well. \hfill\cqd

\nd The next result is an immediate consequence of the proof of  Proposition \ref{raizT_K}. 
\begin{cor}\label{r_0}
Assume $\alpha_1, \alpha_2 \geq 0$ and $0 \leq\gamma_i+\sigma_i <\ell-1$. Then there are an  $r_0=r^{\varepsilon,\delta}_0>0$ and $0 \neq u_1,u_2\in F$, for each $F\in\mathcal{A}$ given, such that  $T_F(u_1,u_2)=(0,0)$ and $\|(u_1,u_2)\|\leq r_0$ with $r_0$ does not depending on $F$. Besides this, $r_0$ does not depend on $\varepsilon>0$ if $\beta_i=0$ and $\alpha_i \in (0,1]$.
\end{cor}
\nd \proof. The first part of Corrolary was just proved by above arguments. To show that $r_0$ does not depend on $\varepsilon>0$, it is necessary just going back to (\ref{511}) and redo the below estimate, say
\begin{eqnarray*}
	\int_\Omega \frac{\hat{a}_i(x)}{(\vert u_i\vert + \varepsilon)^{\alpha_i}}u_i	&\leq &  \int_{\vert u_i\vert \leq 1}\hat{a}_i\vert u_i\vert^{1-\alpha_i} + \int_{\vert u_i\vert > 1}\hat{a}_i\vert u_i\vert^{1-\alpha_i}\nonumber\\ 
	&\leq & \Vert a_i \Vert_1 + \Vert a_i \Vert_{\infty} \Vert u_i \Vert.
\end{eqnarray*}

\nd Now, following similar arguments like those that were used above, we show the existence of one $r_0>0$ independent of $\varepsilon>0$.\hfill\cqd

\nd After these, we are able to  solve the equation $T(u_1,u_2) = 0$, where  the operator $T$ was  given by Proposition \ref{operador-T}. More,  this zero of $T$ will be a solution of \eqref{auxprob}.  To solve $T(u_1,u_2) = 0$, we have inspired in an idea found in the work of Browder \cite{Browder}.  

\begin{lem}\label{raizT}
	For each small $\epsilon > 0$, it there is a $(u^1,u^2)=(u_{\epsilon,\delta}^1,u_{\epsilon,\delta}^2)\in E$, with $u^1,u^2\neq 0$, such that  
	$
	T(u^1,u^2) = 0.
	$
\end{lem}

\nd \proof.  Let $F_0\in\mathcal{A}$  and define 
$$
V_{F_0}=\left\{(u_1,u_2)\in F\times F~ \big|~F\in \mathcal{A},~F_0\subset F,~T_F(u_1,u_2)= 0,~\|(u_1,u_2)\|\leq r_0\right\},
$$
where $r_0>0$ was defined at Corolary \ref{r_0}, and $\mathcal{A}$ was  defined in (\ref{defin}).

\nd By Proposition  \ref{raizT_K} and Corolary  \ref{r_0},  we have that $V_{F_0}\neq \emptyset$ and  $\overline V^\sigma_{F_0} \subset {B_{r_0}}$, where $\overline V^\sigma_{F_0}$ is the weak closure of  $V_{F_0}$ and  $B_{r_0}$ is the closed ball. So, $\overline V_{F_0}^\sigma$ is weakly compact.

\nd {\bf Claim.} The family 
$$\mathcal{B}:=\left\{\overline V_{F}^\sigma~|~F\in \mathcal{A}\right\}
$$ 
\nd has the finite intersection property.

\nd Indeed,  consider the finite family
$
\big\{\overline V_{F_1}^\sigma,\overline V_{F_2}^\sigma,...,\overline V_{F_p}^\sigma\big\}\subset\mathcal{B}
$
and let   $F:=\mbox{span}\{F_1,F_2,...,F_p\}$. So, by the definition of  $V_{F_i}$,  $u_F\in\overline V_{F_i}^\sigma,~i=1,2,...,p$, and so
$$
\bigcap_{i=1}^p\overline V_{F_i}^\sigma\neq\emptyset,
$$
\nd showing the {\bf Claim}.

\nd Since $B_{r_0}$ is weakly compact and $\mathcal{B}$ has the finite intersection property, it follows by  \cite[Thm. 26.9]{munkres} that
$$
W:=\bigcap_{F\in\mathcal{A}}\overline V_{F}^\sigma\neq \emptyset.
$$
\nd Let $(u^1,u^2)=(u_{\epsilon,\delta}^1,u_{\epsilon,\delta}^2)\in W$. Then $T(u^1,u^2)=0$, or equivalently,
\begin{eqnarray}\label{solfrac-+}
	\int_\Omega\phi(|\nabla u^1|)\nabla u^1\nabla\psi dx &=& \int_\Omega\frac{\hat{a}_1(x)\psi}{(|u^1|+\epsilon)^{\alpha_1}(|u^2|+\epsilon)^{\beta_1}}dx\nonumber\\
&+& \int_\Omega \hat{b}_1(x)(u^{1+}+\delta)^{\gamma_1}(u^{2+}+\delta)^{\sigma_1}\psi dx
\end{eqnarray}
\nd and
\begin{eqnarray}\label{solfrac-2+}
	\int_\Omega\phi(|\nabla u^2|)\nabla u^2\nabla\varphi dx &=& \int_\Omega\frac{\hat{a}_2(x)\varphi}{(|u^1|+\epsilon)^{\beta_2}(|u^2|+\epsilon)^{\alpha_2}}dx\nonumber\\
	&+& \int_\Omega \hat{b}_2(x)(u^{1+}+\delta)^{\sigma_2}(u^{2+}+\delta)^{\gamma_2}\varphi dx
\end{eqnarray}
\nd for all $\psi,\varphi\in\w$.
\vskip.1cm

\nd Indeed, take $F_0=\mbox{span}\{\omega_1,\omega
_2,\psi,\varphi,u^1,u^2\}$. So, $ F_0\in \mathcal{A}$ and  by  \cite[Thm. 1.5]{Figueiredo} there is a sequence  $ (u_n^1,u_n^2) \subseteq V_{F_0}$ such that  $(u_n^1,u_n^2)\rightharpoonup (u^1,u^2)$ in $E$. Since  $(u_n^1,u_n^2)\in V_{F_0}$, it folows from its definition that  $\|(u_n^1,u_n^2)\|\leq r_0$ and there exists a $F_n\in \mathcal{A}$, with $F_0\subset F_n$, such that
\begin{eqnarray}\label{u_n}
	\int_\Omega\phi(|\nabla u_n^1|)\nabla u_n^1\nabla \eta_1 dx &=& \int_\Omega\frac{\hat{a}_1(x)\eta_1}{(|u^1_n|+\epsilon)^{\alpha_1}(|u^2_n|+\epsilon)^{\beta_1}}dx\nonumber\\
	&+& \int_\Omega \hat{b}_1(x)(u^{1+}_n+\delta)^{\gamma_1}(u^{2+}_n+\delta)^{\sigma_1}\eta_1 dx,
\end{eqnarray}
and 
\begin{eqnarray}\label{u_n-1}
\int_\Omega\phi(|\nabla u_n^2|)\nabla u_n^2\nabla \eta_2 dx &=& \int_\Omega\frac{\hat{a}_2(x)\eta_2}{(|u^1_n|+\epsilon)^{\beta_2}(|u^2_n|+\epsilon)^{\alpha_2}}dx\nonumber\\
&+& \int_\Omega \hat{b}_2(x)(u^{1+}_n+\delta)^{\sigma_2}(u^{2+}_n+\delta)^{\gamma_2}\eta_2 dx
\end{eqnarray}

\nd for all $\eta_1,\eta_2\in F_n$.
\vskip.2cm

\nd Now, since $\w\stackrel{cpt}\hookrightarrow L_\Phi(\Omega)\hookrightarrow L^\ell(\Omega)$, we have that:
\begin{itemize}
	\item[(1)] $u_n^1\rightarrow u^1$ and $u_n^2\rightarrow u^2$ in $L_\Phi(\Omega)$ and $L^\ell(\Omega)$,
	\item[(2)] $u_n^1\rightarrow u^1$ and $u_n^2\rightarrow u^2$ a.e. in $\Omega$.
\end{itemize}
\nd Set $\eta_i=u_n^i-u^i\in F_n$ for $i=1,2$. So,  it follows from the definition and properties of $\mathcal{A}$ that $\eta_i\in F_n$. Now, by using in $\eta_1$ as a test function in \eqref{u_n}, we get 
\begin{eqnarray}\label{S+}
	\lim\langle (-\Delta_{\Phi})(u_n^1),u_n^1-u^1\rangle&=&\lim\int_\Omega \frac{\hat{a}_1(x)(u_n^1-u^1)}{(|u_n^1|+\epsilon)^{\alpha_1}(|u_n^2|+\epsilon)^{\beta_1}}dx\nonumber\\
	& + & \lim \int_\Omega \hat{b}_1(x)(u^{1+}_n+\delta)^{\gamma_1}(u^{2+}_n+\delta)^{\sigma_1}(u_n^1-u^1) dx \nonumber\\
	&\leq& \lim\int_\Omega \frac{\hat{a}_1(x)}{\epsilon^{\alpha_1+\beta_1}}|u_n^1-u^1|dx.\\
	& + & \lim \int_\Omega \hat{b}_1(x)(|u^{1+}_n|+|u^{2+}_n|+\delta)^{\gamma_1+\sigma_1}|u_n^1-u^1|dx.\nonumber
\end{eqnarray}
\nd Since $\w\stackrel{cpt}\hookrightarrow L^1(\Omega)$, it follows from (1) above, that 
\begin{eqnarray}\label{a-tende0}
	\int_{\Omega}\frac{\hat{a}_1(x)}{\epsilon^{\alpha_1+\beta_1}}|u_n^1-u^1|dx\leq\frac{|\hat{a}_1|_\infty}{\epsilon^{\alpha_1+\beta_1}}\int_{\Omega}|u_n^1-u^1|dx\rightarrow0,
\end{eqnarray}
and by using  $\gamma_1+\sigma_1<\ell-1$,  $\w\hookrightarrow L^\ell(\Omega)$ and (1) again, we obtain
\begin{eqnarray}\label{b-tende0}
\displaystyle\int_\Omega \hat{b}_1(x)(|u^{1+}_n|\!\!\!&+&\!\!\!|u^{2+}_n|+\delta)^{\gamma_1+\sigma_1}|u_n^1-u^1|dx\nonumber\\
	&\leq& |\hat{b}_1|_\infty\left(\int_\Omega(|u^{1+}_n|+|u^{2+}_n|+\delta) ^{\frac{(\gamma_1+\sigma_1)\ell}{\ell-1}}\right)^{\frac{\ell-1}{\ell}}|u_n^1-u^1|_\ell\\
	&\leq& |\hat{b}_1|_\infty\left(\int_\Omega(|u^{1+}_n|+|u^{2+}_n|+\delta) ^\ell\right)^{\frac{\ell-1}{\ell}}|u_n^1-u^1|_\ell.\nonumber
\end{eqnarray}
So, as a consequence of  \eqref{S+}, \eqref{a-tende0} and \eqref{b-tende0}, we have
$$
\lim\langle (-\Delta_{\Phi})(u_n^1),u_n^1-u^1\rangle\leq 0,
$$
\nd which implies that  $u_n^1\rightarrow u^1$ in $\w$, because  $(-\Delta_{\Phi})$ is an operator of the type $(S_+)$ (see \cite[Prop. A.2]{JVMLED}). By a similar argument one shows that  $u_n^2\rightarrow u^2$ in $\w$. 
\smallskip

\nd {\bf Verification of  \eqref{solfrac-+} and \eqref{solfrac-2+}:} Passing to a subsequence, if necessary, we have
\begin{itemize}
	\item[(1)] $\nabla u_n^i\rightarrow \nabla u ^i$ a.e in $\Omega$,
	\item[(2)] there exist $h_i\in L_\Phi(\Omega)$ such that $|\nabla u_n^i|\leq h_i$.
\end{itemize}
\nd By the Young's inequality for N-functions, we have
\begin{eqnarray}
	|\phi(|\nabla u_n^1|)\nabla u_n^1\nabla\varphi|&\leq& \phi(|\nabla u_n^1|)|\nabla u_n^1||\nabla\varphi|\leq \phi(h_1)h_1|\nabla\varphi|\nonumber\\
	&\leq& \widetilde{\Phi}(\phi(h_1)h_1)+\Phi(|\nabla\varphi|)\nonumber\\
	&\leq& \Phi(2h_1)+\Phi(|\nabla\varphi|)\in L^1(\Omega),\nonumber
\end{eqnarray}
\nd where we have used the hypothesis $(\phi_2)$, that is, the function $t\mapsto t\phi(t)$ is increasing for $t\geq 0$.

\nd  Now, applying the Lebesgue's Theorem one finds that
$$
\int_\Omega\phi(|\nabla u_n^1|)\nabla u_n^1\nabla\psi dx \longrightarrow\int_\Omega\phi(|\nabla u^1|)\nabla u^1\nabla\psi dx,
$$	
and arguing as above we also have 
$$
\int_\Omega\phi(|\nabla u_n^2|)\nabla u_n^2\nabla\varphi dx \longrightarrow\int_\Omega\phi(|\nabla u^2|)\nabla u^2\nabla\varphi dx.
$$
\nd Setting $\eta_1=\psi$ and $\eta_2=\varphi$ in \eqref{u_n} and \eqref{u_n-1}, respectively, and passing to the limit, we obtain \eqref{solfrac-+} and \eqref{solfrac-2+}, that is,   $T(u^1,u^2)=(0,0)$. Finally, by similar arguments to those used in Proposition  \ref{raizT_K}, we infer that $u^1,u^2\not\equiv 0$ as well.  \hfill\cqd

\nd To finish the proof of Theorem \ref{solu-aux}, it just remains to show that $u^1,u^2\geq 0$. Taking $-(u^{1})^{-}$  as a test function in Equation \eqref{solfrac-+}, it follows from $(\Phi_3)$ that
\begin{eqnarray}
		\ell\int_\Omega\Phi(|\nabla u^{1-}|)dx&\leq& \int_\Omega\phi(|\nabla u^{1-}|)|\nabla u^{1-}|^2dx\nonumber\\
		& = & -\int_\Omega \frac{\hat{a}_1(x)u^{1-}}{(u^{1-}+\epsilon)^{\alpha_1}(|u^{2}|+\epsilon)^{\beta_1}} dx\nonumber\\
		& - & \int_\Omega \hat{b}_1(x)(u^{1+}+\delta)^{\gamma_1}(u^{2+}+\delta)^{\sigma_1}u^{1-} dx\leq 0,\nonumber
\end{eqnarray}
that is, $u^{1-}\equiv0$. In a similar way one shows that   $u^{2-}\equiv0$. This ends our proof.\hfill\cqd

\section{Proof of Theorem \ref{Teor-prin} (final arguments).}

In this section, we will consider $a_n^i=\max\{a_i,n\}$, and $b_n^i=\max\{b_i,n\}$ for $n \in \mathbb{N}$. So, $a_n^i,b_n^i \in L^{\infty}(\Omega)$ for all $n \in \mathbb{N}$ in accordance to apply Theorem \ref{solu-aux}.
\smallskip

\nd {\it Proof of $(i)$}: In this case, Problem (\ref{auxprob}) read as 
\begin{equation}\label{auxprob-n}
\left\{\
\begin{array}{l}
\displaystyle-\Delta_\Phi u=\frac{a^n_1(x)}{(u+\frac1n)^{\alpha_1}}+b^n_1(x)u^{\gamma_1}v^{\sigma_1},\\ \displaystyle-\Delta_\Phi v=\frac{a^n_2(x)}{(v+\frac1n)^{\alpha_2}}+b^n_2(x)u^{\sigma_2}v^{\gamma_2},~\mbox{in}~\Omega,\\ 
u,v\geq0~\mbox{in}~\Omega,~u=v=0~\mbox{on}~\partial \Omega,
\end{array}
\right.
\end{equation}
\nd where  $\delta=0$,  $\varepsilon=1/n$, $n \in \mathbb{N}$. So, it follows from Theorem \ref{solu-aux} that there exists an $ (u_n,v_n)\in E$, with $u_n,v_n\neq 0$, solution of   (\ref{auxprob-n}).
\smallskip

\nd {\bf Claim.}  $u_n+ {1}/{n}\geq C d(x) $  and $v_n+ {1}/{n}\geq C d(x) $ for some $C>0$. 

\nd Indeed,  by taking $\hat{b}_i=0$ and applying Theorem \ref{solu-aux} with $\varepsilon=1$, we obtain a $ w_i \in \w$ solution of the problem 
\begin{equation}\label{aux-1}
	\left\{
	\begin{array}{ll}
		-\Delta_\Phi w=\displaystyle\frac{a^1_i(x)}{(w+1)^{\alpha_i}}~\mbox{in}~\Omega,\\ 
		w\geq 0~\mbox{in}~ \Omega;~ w=0~\mbox{on}~\partial \Omega.
	\end{array}
	\right.
\end{equation}

\nd Besides this, it follows from \cite[Lemma 3.3]{Fukagai}, that $w_i\in C^{1,\beta_i}(\overline{\Omega})$, for some $0<\beta_i<1$, and from  \cite[Prop. 5.2]{JVMLED}, it follows that $w_i>0$.

\nd Since the solution $(u_n,v_n)$ of \eqref{auxprob-n} satisfies
\begin{equation*}\label{super}
-\Delta_\Phi u_n\geq \frac{a^n_1(x)}{(u_n+1/n)^{\alpha_1}}\geq \frac{a^1_1(x)}{(u_n+1)^{\alpha_1}}~\mbox{in}~\Omega,
\end{equation*}
it follows that $w_1$ and $u_n$ are respectively sub e supersolutions of \eqref{aux-1} with $i=1$ such that
$$
0\leq \frac{w_1}{u_n+\frac{1}{n}}\leq  nw_1\in L^\infty(\Omega).
$$

\nd So, by Theorem \ref{princ-compar}, we obtain

\begin{equation}\label{u-posi}
 u_n+\frac{1}{n}\geq w_1>0~\mbox{a.e. in}~\Omega.
\end{equation}
Finally, by a classical argument, we can show that $w_1 \geq C_1d$, for some $C_1>0$. In a similar way, we have 
\begin{equation}\label{v-posi}
	v_n(x)+\frac 1 n\geq w_2(x)\geq C_2 d(x)~a.e.~ \mbox{in}~ \Omega,
\end{equation}
as well. This ends the proof of Claim.

\nd Now, by using $u_n\in \w$ as a test function in the first equation in (\ref{auxprob}), the hypothesis $(\phi_3)$ ($(\phi_3)$ implies $(\phi_3)^{\prime}$, see Remark \ref{8.1}), the arguments as in (\ref{551}), and (\ref{u-posi}), we obtain
\begin{eqnarray}\label{bound}
\ell \zeta_0(\Vert \nabla u_n \Vert_{\Phi}) &\leq& \ell \int_{\Omega} \Phi(\vert \nabla u_n \vert) \leq \int_{\Omega} \phi(\vert \nabla u_n \vert) \vert \nabla u_n \vert^2\nonumber\\
	&\leq &\int_{\Omega}\frac{a_1(x)}{(u_n+ 1/n)^{\alpha_1}}u_ndx + \int_{\Omega}b_1(x) u_n^{\gamma_1}v_n^{\sigma_1} u_ndx\\
	&=&\int_{\Omega}\frac{a}{d^{\alpha_1}}u_ndx+ \int_{\Omega}b_1 (v_n+ u_n)^{\gamma_1+\sigma_1+1} dx\nonumber\\
	&\leq& \Vert a/d^{\alpha_1} \Vert_{\tilde{\Psi}}\|u_n\|_\Psi+C_2' \Vert b_1 \Vert_{q_1}\Vert (u_1,u_2) \Vert^{\gamma_1+\sigma_1 + 1}
	\nonumber\\
	&\leq& C_1\|u_n\|+C_2\Vert (u_1,u_2) \Vert^{\gamma_1+\sigma_1 + 1},\nonumber
	\end{eqnarray}
and, in a similar way, we obtain
\begin{eqnarray}\label{bound1}
\ell \zeta_0(\Vert \nabla v_n \Vert_{\Phi}) 
	\leq D_1\|v_n\|+D_2\Vert (u_1,u_2) \Vert^{\gamma_2+\sigma_2 + 1}.
	\end{eqnarray}

\nd So, by using either (\ref{bound}) or (\ref{bound1}), we can show that $(u_n,v_n) \subset E$ is bounded if either $\Vert v_n \Vert \geq 1$ and $\Vert u_n \Vert \leq 1$ or $\Vert v_n \Vert \leq 1$ and $\Vert u_n \Vert \geq 1$. Now, assume that $\Vert u_n \Vert , \Vert u_n \Vert  \geq 1$. So, by summing (\ref{bound}) and (\ref{bound1}), it follows from Lemma \ref{lema_naru}, that 
$$
\Vert u_n \Vert^{\ell} + \Vert v_n \Vert^{\ell} 
	\leq D_1'\|(u_n,v_n)\|+D_2'\Vert (u_1,u_2) \Vert^{\gamma_1+\sigma_1 + 1}+ D_3'\Vert (u_1,u_2) \Vert^{\gamma_2+\sigma_2 + 1}.
$$
for some $D_i'>0$. Since, $\gamma_i+\sigma_i < \ell -1$, we have that $(u_n,v_n)\subset E$ is bounded as well. Passing to a subsequence if necessary, we find that
\begin{itemize}
	\item[(1)] $u_n\rightharpoonup u$ and $v_n\rightharpoonup v$ in $\w$;
	\item[(2)] $u_n\rightarrow u$ and $v_n\rightharpoonup v$ in $L_\Phi(\Omega)$;
	\item[(3)] $u_n\rightarrow u$ e $v_n\rightharpoonup v$ a.e. in $\Omega$;
	\item[(4)] there exist $\theta_i\in L_\Phi(\Omega)$ such that  $0<C_1 d\leq u_n + 1/n\leq \theta_1$ and $0<C_2 d\leq v_n+1/n\leq \theta_2$,
\end{itemize}
that is, by using (3), \eqref{u-posi} and \eqref{v-posi} we have that  $u,v\geq Cd$ a.e. in $\Omega$ for some $C>0$.

\nd Next we will show that $(u_n,v_n)\rightarrow(u,v)$ in $E$. Indeed, by using  (4) and $a/d^{\alpha_1} \in L_{\tilde{\Psi}}$, we obtain
\begin{eqnarray}
	\int_{\Omega}\frac{a^n_1(u_n-u)}{(u_n+\frac{1}{n})^{\alpha_1}}&\leq & \int_{\Omega}\frac{a_1}{d^{\alpha_1}}\vert u_n -u\vert\leq 2\Vert{a_1}/{d^{\alpha_1}}\Vert_{\tilde{\Phi}} \Vert u_n - u \Vert_{\Phi}\nonumber  
\end{eqnarray}
\nd and
$$|b_1(x)u_n^{\gamma_1}v_n^{\sigma_1}(u_n-u)|\leq 2b_1(x)\max\{\theta_1,\theta_2\}^{\gamma_1+\sigma_1+1}\in L^{1}(\Omega),$$
because $b_1 \in L^{q_1}(\Omega)$ and $L^{\Phi}(\Omega)\hookrightarrow L^{\gamma_1+\sigma_1 + 1}(\Omega)$. So, by applying Fatou's Lemma, it follows from \eqref{auxprob-n} that
$$
\limsup\langle (-\Delta_{\Phi})u_n,u_n-u\rangle\leq 0,
$$
that is, $u_n\rightarrow u$ in $\w$. Similarly one shows that  $v_n\rightarrow v$ em $\w$.  Passing to the limit in \eqref{auxprob-n} we infer that $(u,v)$ is a weak solution of \eqref{prob}. These end the proof of $(i)$. 
\smallskip

\nd About Remark \ref{remark1}, it is necessary just to observe that $r_0$ given  by  Corrolary \ref{r_0} does not depend on $n$, that is, $(u_n,v_n) \subset E$ is already bounded. So, the remaining arguments are the same.
\medskip

\nd {\it Proof of $(ii)$}:  In this case, the problem (\ref{auxprob}) reduces to 
\begin{equation}\label{prob-aux-2b}
\left\{\
\begin{array}{l}
\displaystyle-\Delta_\Phi u=\frac{a^n_1(x)}{(u+ 1/n)^{\alpha_1}(v+1/n)^{\beta_1}}+b^n_1(x)(u+1/n)^{\gamma_1}~\mbox{in}~\Omega,\\
\displaystyle-\Delta_\Phi v=\frac{a^n_2(x)}{(u+1/n)^{\beta_2}(v+1/n)^{\alpha_2}}+b^n_2(x)(v+1/n)^{\gamma_2}~\mbox{in}~\Omega,\\ 
u,v \geq 0~\mbox{in}~\Omega,~u=v=0~\mbox{on}~\partial \Omega,
\end{array}
\right.
\end{equation}
by taking $\varepsilon=\delta=1/n$ with $n\in \mathbb{N}$. So, it follows from  Theorem \ref{solu-aux} that there exists an $ (u_n,v_n)\in E$, with $u_n,v_n\neq 0$, solution of   (\ref{prob-aux-2b}).

\nd Besides this, with similar arguments as those used in Case $(i)$, we able to show that
$$u_n+{1}/{n}\geq Cd,~  v_n+1/ n \geq Cd,~\mbox{for some}~C>0 .$$
In this case, we redo the above arguments using $w_i \in C^{1,\tau_i}(\overline{\Omega})$,  for some $0<\tau_i<1$, as solution of the problem
\begin{equation*}\label{solu-aux-2}
\left\{\
\begin{array}{l}
\displaystyle-\Delta_\Phi w=b^1_i(x)w^{\gamma_i}~ \mbox{in}~\Omega, \\
w>0~\mbox{in}~\Omega,~w=0~\mbox{on}~\partial \Omega,
\end{array}
\right.
\end{equation*}
in what the existence result is given by Theorem \ref{solu-aux} with $\hat{a_i}=0$, and $\delta=0$. The regularity is guaranteed by  \cite[Corolary 3.1]{fang}, and the positivity is given by  \cite[Prop. 5.2]{JVMLED} again.

\nd After this, in a similar way to those that we have done in case $(i)$, we are able to show that $(u_n,v_n) \subset E$ is bounded as well. This ends the proof of Theorem \ref{Teor-prin} - $(ii)$. 
\vskip.2cm

\nd {\it Proof of $(iii)$}:   Under these  conditions,  the system  (\ref{auxprob}) becomes
\begin{equation}\label{prob-aux-3b1}
\left\{\
\begin{array}{l}
\displaystyle-\Delta_\Phi u=\frac{a^n_1(x)}{(v+1/n)^{\beta_1}}+b^n_1(x)(v+1/n)^{\sigma_1},\\

\displaystyle-\Delta_\Phi v=\frac{a^n_2(x)}{(u+1/n)^{\beta_2}}+b^n_2(x)(u+1/n)^{\sigma_2},~\mbox{in}~\Omega,\\ 
u,v\geq0~\mbox{in}~\Omega,~u=v=0~\mbox{on}~\partial \Omega.
\end{array}
\right.
\end{equation}
\nd with $\varepsilon=\delta=1/n$ and $n \in \mathbb{N}$.

\nd So, it follows from Theorem \ref{solu-aux} that there exists 
a week solution $(u_n,v_n)\in E$ to Problem \eqref{prob-aux-3b1}.

\nd {\bf Claim.} there exists a constant $C>0$ such that $u_n\geq C d~ \mbox{and}~  v_n\geq C d.$

\nd Indeed, let us set $h_i(x):=\min\{a_i(x),b_i(x)\}$, $x \in \Omega$, and 
$$
g_i:=\min_{t\geq 0} \left[\frac{1}{t^{\beta_i}+1}+t^{\sigma_i}\right].
$$
So, we can  infer from our hypotheses that 
$$h_i(x)\geq 0,~  h_i\not\equiv 0~\mbox{and}~g_i>0.$$

\nd Now, it follows from \cite[Lemma 3.4]{Fukagai} that there is a positive solution  $w_i\in\w$ of
\begin{equation}\label{auxprob-2-3}
\left\{\
\begin{array}{l}
\displaystyle-\Delta_\Phi w= h_i(x)g_i~\mbox{in}~\Omega,\\ 
w\geq 0~\mbox{in}~\Omega,~w=0~\mbox{on}~\partial \Omega.
\end{array}
\right.
\end{equation}
Again, it follows from \cite[Corolary [3.1]{fang} that $w_i\in C^{1,\tau_i}(\overline{\Omega})$, for some $0<\tau_i<1$, and by \cite[Prop. 5.2]{JVMLED}, we have $w_i>0$, thais is, there exists a $C_i>0$ such that $w_i \geq C_i d$.

\nd So, it follows from (\ref{prob-aux-3b1}) and (\ref{auxprob-2-3}), that 
$$
\left\{\
\begin{array}{l}
\displaystyle-\Delta_\Phi u_n\geq h_1(x)g_1=-\Delta_\Phi w_1~\mbox{in}~\Omega,\\ 
\displaystyle-\Delta_\Phi v_n\geq h_2(x)g_2=-\Delta_\Phi w_2~\mbox{in}~\Omega,\\
u_n=v_n=w_1=w_2=0~\mbox{on}~\partial \Omega,
\end{array}
\right.
$$
that implies, by applying  \cite[Lemma 4.1]{fang},  that $u_n\geq w_1\geq C_1 d$ and $v_n\geq w_2\geq C_2 d$, for some $C_1,C_2>0$.

\nd Finally, arguing as in the proof of  Theorem \ref{Teor-prin}-$(i)$, we are able to show that  $(u_n,v_n) \subset E$ is bounded as well. These end the proof of Theorem \ref{Teor-prin} - $(iii)$ and the proof of Theorem \ref{Teor-prin}.\hfill\cqd

\section{Appendix}

\nd  In this section we present for, the reader's convenience, several results used in this paper. We begin referring the reader  to  $\cite{A,Rao1}$ regarding to Orlicz-Sobolev spaces.  The usual norm on $L_{\Phi}(\Omega)$ is 
\[
\|u\|_\Phi=\inf\left\{\lambda>0~|~\int_\Omega \Phi\left(\frac{u(x)}{\lambda}\right) dx \leq 1\right\},
\]
called as Luxemburg norm, and the   norm on $ W^{1, \Phi}(\Omega)$ is
\[
\displaystyle \|u\|_{1,\Phi}=\|u\|_\Phi+\sum_{i=1}^N\left\|\frac{\partial u}{\partial x_i}\right\|_\Phi,
\]
known as Orlicz-Sobolev norm, while $\w$ stands for the closure of $C_0^{\infty}(\Omega)$ with respect to the  norm in $W^{1,\Phi}(\Omega)$. 
\nd We also recall that, under hypotheses $(\Phi_1)-(\Phi_3)$, the functions  $\Phi$ and $\widetilde{\Phi}$  are  N-functions  satisfying  the $\Delta_2$-condition (see \cite[pg 22]{Rao1}). As consequence of these, we have  $L_{\Phi}(\Omega)$  and $W^{1,\Phi}(\Omega)$  are separable, reflexive, and Banach spaces.

\begin{rmk}
\label{8.1}
	It is well known that $(\phi_3)$ implies that 
	\begin{itemize}
		\item[$(\phi_3)^\prime$] $\displaystyle\ell\leq \frac{\phi(t)t^2}{\Phi(t)}\leq m,~t>0,$
	\end{itemize}
	is verified.
	Furthermore, these hypotheses imply that $\Phi,\widetilde\Phi\in\Delta_2$.
\end{rmk}

\nd The inequality 
\[
\int_\Omega\Phi(u)dx\leq \int_\Omega\Phi(2\overline{d}|\nabla u|)dx~\mbox{for all}~ u \in \w,
\]
\nd where $\overline{d}>0$ is the diameter of $\Omega$, is known as  Poincar\'e's inequality (see e.g.  \cite{gossez-Czech}), and as a consequence of it, we have  
\[
\|u\|_\Phi\leq 2\overline{d}\|\nabla u\|_\Phi~\mbox{for all}~ u \in \w,
\]
that is,  $\|u\| :=\|\nabla u\|_\Phi$ defines an equivalent norm to the $\|.\|_{1,\Phi}$  in $\w$. 

\nd Now, let $\Phi_*$ be the inverse of the function
$$
t\in(0,\infty)\mapsto\int_0^t\frac{\Phi^{-1}(s)}{s^{\frac{N+1}{N}}}ds,
$$
\nd which is extended to $\mathbb{R}$ by  $\Phi_*(t)=\Phi_*(-t)$ for  $t\leq 0.$
We say that an N-function $\Psi$ grows essentially more slowly than $\Phi_*$, we denote this by $\Psi<<\Phi_*$, if
$$
\lim_{t\rightarrow \infty}\frac{\Psi(\lambda t)}{\Phi_*(t)}=0~\mbox{for all}~\lambda >0.
$$

\nd The imbeddings below (see \cite{A}) was  used in this paper
$$
\displaystyle \w \stackrel{\tiny cpt}\hookrightarrow L_\Psi(\Omega)~\mbox{if}~\Psi<<\Phi_*,
$$
and in particular, as $\Phi<<\Phi_*$ (see \cite[Lemma 4.14]{Gz1}), we have
$$
\w \stackrel{\tiny{cpt}} \hookrightarrow L_\Phi(\Omega), ~ \mbox{and}~ W_0^{1,\Phi}(\Omega) \stackrel{\mbox{\tiny cont}}{\hookrightarrow} L_{\Phi_*}(\Omega).
$$

\nd It is worth to note that if just $(\phi_1)-(\phi_2)$, and $(\phi_3)^\prime$ are verified, then
$$L_\Phi(\Omega)\stackrel{\mbox{\tiny cont}}\hookrightarrow L^\ell(\Omega)$$
is true, see \cite[Lemma. D.2]{clement}.

\nd Below, we state some Lemmas whose proofs can be find in \cite{Fuk_1}.

\begin{lem}\label{lema_naru}
	Assume  $\phi$ satisfies  $(\phi_1)-(\phi_3)$.
	Set
	$$
	\zeta_0(t)=\min\{t^\ell,t^m\},~\mbox{and}~ \zeta_1(t)=\max\{t^\ell,t^m\},~ t\geq 0.
	$$
	\nd Then  $\Phi$ satisfies
	$$
	\zeta_0(t)\Phi(\rho)\leq\Phi(\rho t)\leq \zeta_1(t)\Phi(\rho),~ \rho, t> 0,
	$$
	$$
	\zeta_0(\|u\|_{\Phi})\leq\int_\Omega\Phi(u)dx\leq \zeta_1(\|u\|_{\Phi}),~ u\in L_{\Phi}(\Omega).
	$$
\end{lem}
\begin{lem}\label{lema_naru_*}
	Assume that  $\phi$ satisfies $(\phi_1)-(\phi_3)$.  Set
	$$
	\zeta_2(t)=\min\{t^{\widetilde\ell},t^{\widetilde m}\},~\mbox{and}~ \zeta_3(t)=\max\{t^{\widetilde\ell},t^{\widetilde m}\},~  t\geq 0,
	$$
	\nd where $1<\ell,m<N$, $\widetilde m = {m}/{(m-1)}$, and $\widetilde\ell = {\ell }{(\ell-1)}$.  Then
	$$
	\widetilde\ell\leq\frac{t^2\widetilde\Phi'(t)}{\widetilde\Phi(t)}\leq \widetilde m,~t>0,~\zeta_2(t)\widetilde\Phi(\rho)\leq\widetilde\Phi(\rho t)\leq \zeta_3(t)\widetilde\Phi(\rho),~ \rho, t> 0,
	$$
and
	$$
	\zeta_2(\|u\|_{\widetilde\Phi})\leq\int_\Omega\widetilde\Phi(u)dx\leq \zeta_3(\|u\|_{\widetilde\Phi}),~ u\in L_{\widetilde\Phi}(\Omega).
	$$
\end{lem}

\begin{lem}\label{lema_Phi}
Let  $\Phi$ be an $N$-function satisfying  $\Delta_2$ condition.  Let  $(u_n)\subset L_\Phi(\Omega)$ be a sequence such that $u_n\rightarrow u$ in  $L_\Phi(\Omega)$. Then there is a subsequence $(u_{n_k})\subseteq(u_n)$ such that:
	\begin{description}
		\item{\rm{(i)}} $u_{n_k}(x)\rightarrow u(x) $ a.e. $x\in\Omega$,
		\item{\rm{(ii)}} there is  an $h\in L_\Phi(\Omega)$ such that 
		$|u_{n_k}|\leq h~\mbox{a.e. in}~ \Omega.$
	\end{description}
\end{lem}
\proof ~ (Sketch): Since $L_\Phi(\Omega)\hookrightarrow L^1(\Omega)$, (see \cite{A}), passing to a subsequence if necessary, we have
$u_n\rightarrow u~\mbox{a.e.}~\mbox{in}~\Omega$.
\nd Moreover, since
$$\int_\Omega\Phi(u_n-u)dx\rightarrow 0,$$ 
 there exists an $\widetilde{h}\in L^1(\Omega)$ such that $\Phi(u_n-u)\leq\widetilde{h}~\mbox{a.e.}~\mbox{in}~\Omega$.
\vskip.1cm

\nd Now, by using that $\Phi$ is convex, increasing  and satisfies $\Delta_2$ condition, we have
$$
\begin{array}{lll}
\Phi(|u_n|) & \leq & \displaystyle C\Phi\left(\frac{|u_n-u|+|u|}{2}\right) \leq \displaystyle\frac{C}{2}\left[\Phi(|u_n-u|)+\Phi(|u|)\right] \leq\displaystyle\frac{C}{2}[\widetilde{h}+\Phi(|u|)].
\end{array}
$$
\nd Defining $h=\Phi^{-1}\big({C}/{2}(\widetilde{h}+\Phi(|u|))\big)$, it follows from $\widetilde{h}\in L^1(\Omega)$ and $\Phi(|u|)\in L^1(\Omega)$, that
$$
\begin{array}{lll}
\displaystyle\int_\Omega\Phi(h)dx & = & \displaystyle\int_\Omega\Phi\left(\Phi^{-1}\left(\frac{K}{2}(\widetilde{h}+\Phi(|u|))\right)\right)dx \\
& = & \displaystyle \int_\Omega\left(\frac{K}{2}(\widetilde{h}+\Phi(|u|))\right)dx <  \infty,
\end{array}
$$
\nd showing that $h\in L_\Phi(\Omega)$. \hfill\cqd

\end{document}